%% file: main.tex
\documentclass[a4paper,parskip=half,11pt,oneside]{scrartcl}
\usepackage{amsmath}
\usepackage[utf8]{inputenc}
\usepackage[OT2,T1]{fontenc}
\usepackage[russian,english]{babel}
\usepackage[all]{xy}
\usepackage{mathrsfs}
\usepackage{amssymb}
\usepackage{booktabs}
\usepackage{palatino}
\usepackage{enumitem}
\usepackage{eulervm}
\usepackage{tikz}
\usepackage[amsmath,thmmarks]{ntheorem}
\usepackage{url}
\usepackage{fancyhdr}
\usepackage{setspace}
\usepackage{hyperref}
\usepackage[capitalize]{cleveref}
\usepackage[bindingoffset=0mm,inner=28mm,outer=28mm,bottom=46mm]{geometry}

\theoremheaderfont{\bf\sffamily}
\newtheorem{Theorem}{Theorem}[section]
\newtheorem{Lemma}[Theorem]{Lemma}
\newtheorem{Proposition}[Theorem]{Proposition}

\theorembodyfont{\normalfont}
\newtheorem{Definition}[Theorem]{Definition}
\newtheorem{Remark}[Theorem]{Remark}
\newtheorem{Question}[Theorem]{Question}
\theoremsymbol{\ensuremath{\diamond}}
\newtheorem{Example}[Theorem]{Example}
\theoremstyle{nonumberplain}
\theoremheaderfont{\sc}
\theoremseparator{:}
\theoremsymbol{\ensuremath{\square}}
\newtheorem{Proof}{Proof}

\crefname{Theorem}{Theorem}{Theorems}
\crefname{Lemma}{Lemma}{Lemmata}
\crefname{Example}{Example}{Examples}
\crefname{Proposition}{Proposition}{Proposition}
\crefname{Remark}{Remark}{Remark}

\newcommand{\fra}{\mathfrak{a}}
\newcommand{\frB}{\mathfrak{B}}

\newcommand{\CH}{A^\ast}
\newcommand{\bbP}{\mathbb{P}}
\newcommand{\bbQ}{\mathbb{Q}}
\newcommand{\calG}{\mathcal{G}}
\newcommand{\scrO}{\mathscr{O}}
\newcommand{\scrU}{\mathscr{U}}
\newcommand{\F}{\mathcal{F}}
\newcommand{\M}{\mathcal{M}}

\newcommand{\Gr}{\calG}
\DeclareMathOperator{\PGL}{PGL}
\DeclareMathOperator{\Spec}{Spec}

\renewcommand{\phi}{\varphi}

\setenumerate{label=(\alph*)}
\newlist{enumstep}{enumerate}{10}
\setlist[enumstep]{label=\emph{\arabic*.}}

\begin{document}

\author{Mathias Häbich}\date{}
\title{Mustafin degenerations}

\maketitle

\begin{center}
  \begin{minipage}[h]{.7\linewidth}
    {\bf Abstract.}
    A \emph{Mustafin degeneration} is a degeneration of a flag variety induced by a configuration of vertices in the Bruhat-Tits building of the projective linear group over a field with a non-archimedean discrete valuation. In the case where the flag type is projective space, Mustafin degenerations have been studied previously.  We generalize the construction to the case of arbitrary flag types and
study the behavior of the components of the special fiber under a natural projection morphism between Mustafin degenerations, which arises from the inclusion of flags. 
  \end{minipage}
\end{center}

\input{intro}
\input{basics}
\input{combinatorics}

\bibliographystyle{beta}
\bibliography{diss}

\bigskip

\textbf{Author's address:}

Mathias Häbich, Institut für Mathematik, Goethe-Universität, 60325 Frankfurt am Main, Germany, \url{haebich@math.uni-frankfurt.de}

\end{document}

%% file: intro.tex
\section{Introduction}

In this paper, we discuss certain degenerations of flag varieties, which we call \emph{Mustafin degenerations}.
By definition, a Mustafin degeneration is induced by a vertex configuration in the Bruhat-Tits building of the projective linear group over a discretely valued field $K$.
 Since these vertices can be described by matrices with entries in $K$, our degenerations can also be thought of as being induced by sets of such matrices, and thus also arise naturally from the perspective of linear algebra. In addition, this point of view gives rise to explicit equations for Mustafin degenerations.

We give a brief account of the definition. Details can be found at the beginning of \cref{sec:definitions}. Fix a vector space $V$ of finite dimension over a discretely valued field $K$ with valuation ring $R$. A vertex in the Bruhat-Tits building is represented by a lattice $L$ in $V$, to which we can assign the flag variety $\F(L)$ parametrizing flags of a specified type in $L$. Given a finite set $\Gamma$ of vertices, represented by lattices $L_1, \ldots, L_n$, we define the Mustafin degeneration $\M_\F(\Gamma)$ as the join of the flag varieties $\F(L_1), \ldots, \F(L_n)$ along their common generic fiber $\F(V)$.

These schemes $\M_\F(\Gamma)$ are natural generalizations of similarly constructed degenerations of projective space, called \emph{Mustafin varieties} in \cite{mustafinvarieties}, which were introduced by Mumford in his influential work \cite{mumford} on the uniformization of curves, and later generalized by Mustafin \cite{mustafin} to higher dimensions.

 In \cref{thm:basic-properties}, we show that Mustafin degenerations are integral schemes, flat and projective over $R$, with connected special fiber. In \cref{rem:reduced}, we give a sufficient condition for the special fiber to be reduced, and in \cref{lem:number-of-components} we prove an upper bound for the number of components of the special fiber in certain cases.

In \cite{mustafinvarieties}, we introduced a distinction of the irreducible components of the special fiber of a Mustafin variety into primary and secondary components, where primary components correspond to vertices in $\Gamma$ and secondary components correspond to vertices in the convex hull of $\Gamma$ that are not contained in $\Gamma$.
However, for arbitrary flag types, the situation is different, since there may be components that are neither primary nor secondary. We show how such so-called tertiary components appear by studying the behavior of the different kinds of components under a natural projection morphism between Mustafin degenerations, which arises from the inclusion of flags.

This work is organized as follows. In \cref{sec:definitions}, we present our framework and give the definition of a Mustafin degeneration, as well as the basic geometric properties. 
\cref{sec:combinatorics} is devoted to the study of the components of the special fiber of a Mustafin degeneration.


%% file: basics.tex
\section{Definition and First Properties}\label{sec:definitions}

First, we review the definition of the Bruhat-Tits building of the projective linear group over a non-archimedean field. Let $R$ be a discrete valuation ring, $K$ its field of fractions, and $V$ a $K$-vector space of dimension $d$. A \emph{lattice} is a free $R$-submodule of $V$ of rank $d$. Two lattices $L'$ and $L$ are called \emph{homothetic} if $L' = cL$ for some $c \in K^\times$. Homothety is an equivalence relation, and an equivalence class $[L]$ is called a \emph{vertex}. The set of all vertices is denoted by $\frB_d^0$. Two vertices $[L_1]$ and $[L_2]$ are called \emph{adjacent} if there are representatives $L_1' \in [L_1]$ and $L_2' \in [L_2]$ such that $\pi L_1' \subsetneq L_2' \subsetneq L_1'$. This is a symmetric relation and thus defines an undirected graph on the vertex set $\frB_d^0$. The clique complex of this graph, where each set of pairwise adjacent vertices forms a simplex, is a simplicial complex denoted by $\frB_d$. It is called the \emph{Bruhat-Tits building} of $\PGL(V)$. See \cite[Section~6.9]{abramenko-brown} for details of the construction and some properties of $\frB_d$.

Let $\F = (k_1 < k_2 < \cdots < k_r)$ be a tuple of ascending integers satisfying $1 \leq k_i \leq d-1$. We call this a \emph{flag type}.
To a vertex $[L] \in \frB_d^0$, we assign the flag variety $\F(L)$, which is defined as the $R$-scheme representing the functor which sends a scheme $h \colon T \to \Spec R$ to the set
\[\{h^*(\widetilde L) \supset \scrU_1 \supset \scrU_2 \supset \dots \supset \scrU_r  : h^*(\widetilde L)/\scrU_i \text{ is a locally free $\scrO_T$-module of rank $k_i$}\}\text{.}\]
$\F(L)$ is an integral scheme, projective and smooth over $R$. If $\F = (k)$, $\F(L)$ is the Grassmannian $\Gr_k(L)$ parametrizing quotients of dimension $k$, and in particular,  we recover the projective space $\bbP(L) \cong \bbP_R^{d-1}$ (parametrizing hyperplanes) for $k=1$ and the dual projective space $\bbP^\vee(L)$ (parametrizing lines) for $k=d-1$. In these cases, we also write $\F = \Gr_k$, $\F = \bbP$ and $\F = \bbP^\vee$ instead of $\F = (k)$, $\F = (1)$, and $\F = (d-1)$, respectively.

It will become important in \cref{sec:combinatorics} to observe that an inclusion
\[\F' = (k_{i_1} < k_{i_2} < \cdots < k_{i_{r'}}) \subseteq (k_1 < k_2 < \cdots < k_r) = \F\]
of flag types induces a surjective morphism $\F(L) \twoheadrightarrow \F'(L)$, given on $T$-valued points by sending a flag $\{\scrU_1 \supset \scrU_2 \supset \dots \supset \scrU_r\}$ as above to the subflag $\{\scrU_{i_1} \supset \scrU_{i_2} \supset \dots \supset \scrU_{i_{r'}}\}$.

Since, by definition, any two lattices  $L$ and $L'$  are isomorphic as $R$-modules, the associated flag varieties are isomorphic as $R$-schemes. Note, however, that the isomorphism is not canonical, since it depends on a choice of basis for each lattice. On the other hand, their generic fibers $\F(L) \times_R K$ and $\F(L') \times_R K$ are canonically isomorphic to $\F(V)$, the $K$-scheme representing the functor which sends $h \colon T \to \Spec K$ to 
\[\{h^*(\widetilde V) \supset \scrU_1 \supset \scrU_2 \supset \dots \supset \scrU_r  : h^*(\widetilde V)/\scrU_i \text{ is a locally free $\scrO_T$-module of rank $k_i$}\}\text{.}\]

This allows us to make the following definition:
\begin{Definition}\label{def:mustafin-degeneration}
  Given a flag type $\F$ as above and a finite subset $\Gamma = \{[L_1], \ldots, [L_n]\} \subset \frB_d^0$, we define the \emph{Mustafin degeneration} $\M_\F(\Gamma)$ to be the join of the  schemes $\F(L_1), \ldots, \F(L_n)$, which is the scheme-theoretic image (i.e., the closure of the image with the induced reduced subscheme structure) of the map
\begin{equation}\label{eqn:def-mustafin}
  \F(V) \longrightarrow \F(L_1) \times_R \cdots \times_R \F(L_n)\text{.}
\end{equation}
\end{Definition}

We can also describe $\M_\F(\Gamma)$ by giving equations. Recall that the flag variety $\F(L)$ with $\F = (k_1 < \cdots < k_r)$ can be embedded into the product of projective spaces $\prod_{t=1}^r \bbP^{\binom{d}{k_t}-1}_R$ by the Plücker embedding. The equations cutting out the image of this embedding can be found in \cite[Equations~(1) and~(3)]{youngtableaux}. If, for $1 \leq j \leq n$, $e_1^{(j)}, \ldots, e_d^{(j)}$ is a basis for $L_j$, then the respective multihomogeneous coordinates on $\prod_{t=1}^r \bbP^{\binom{d}{k_t}-1}_R$ are
\[ (p_{i_1, \ldots, i_{k_t}}^{(j)} : 1 \leq i_1 < \cdots < i_{k_t} \leq d)_{t=1,\ldots,r}\text{,} \]
 where $p_{i_1, \ldots, i_{k_t}}^{(j)} = e_{i_1}^{(j)} \wedge \dots \wedge e_{i_{k_t}}^{(j)}$. Let $A^{(j)}_t$ be the matrix satisfying $A^{(j)}_tp_{i_1, \ldots, i_{k_t}}^{(1)} = p_{i_1, \ldots, i_{k_t}}^{(j)}$. $\mathbold{A}^{(j)} := (A_1^{(j)}, \ldots, A_r^{(j)})$ is an isomorphism $\F(L_1) \to \F(L_j)$ inducing an automorphism $\F(V) \to \F(V)$. The diagram
\[ \xymatrix@C+5cm{
  \F(V) \ar[r]^{((\mathbold{A}^{(1)})^{(-1)}, \ldots, (\mathbold{A}^{(n)})^{(-1)}) \circ \Delta} \ar[d] & \F(V)^n \ar[d] \\
    \prod_{[L] \in \Gamma} \F(L) \ar[r]^{((\mathbold{A}^{(1)})^{(-1)}, \ldots, (\mathbold{A}^{(n)})^{(-1)})} & \F(L_1)^n
} \]
commutes, where $\Delta \colon \F(V) \to \F(V)^n$ is the diagonal map and the vertical arrow on the left is the map from \cref{eqn:def-mustafin}; hence, the Mustafin degeneration is cut out in $\prod_{[L] \in \Gamma}\prod_{t=1}^r \bbP^{\binom{d}{k_t}-1}_R$ by the ideal $\fra \cap R[\ldots, p_{i_1, \ldots, i_{k_t}}^{(j)}, \ldots]$, where $\fra$ is
the ideal generated over $K$ by all $2 \times 2$-minors of the matrices
\[\begin{pmatrix}
  \vdots & & \vdots \\
  p_{i_1, \ldots, i_{k_t}}^{(1)} & \cdots & (A^{(j)}_t)^{-1} p_{i_1, \ldots, i_{k_t}}^{(j)} & \cdots \\
  \vdots & & \vdots
\end{pmatrix} \qquad (t = 1, \ldots, r)\text{,}\]
together with the equations cutting out the product of flag varieties $\prod_{[L] \in \Gamma} \F(L)$.

Mustafin degenerations have the following geometric properties:

\begin{Theorem}\label{thm:basic-properties}
  For a finite subset $\Gamma$ of $\frB_d^0$, any Mustafin degeneration $\M_\F(\Gamma)$ is an integral scheme which is flat and projective over $R$. Its generic fiber is equal to $\F(V)$, and its special fiber $\M_\F(\Gamma)_s$ is connected and equidimensional of dimension $\dim \F(V)$.
\end{Theorem}

\begin{Proof}
   By construction, the Mustafin degeneration $\M_\F(\Gamma)$ is a reduced, irreducible, projective scheme over $R$ with generic fiber $\F(V)$. Since $R$ is a discrete valuation ring, flatness follows from the fact that $\M_\F(\Gamma)$ is reduced with non-empty generic fiber by \cite[Proposition~4.3.9]{kinglouis}.

In order to show that the special fiber is connected, we apply Zariski's Connectedness Principle \cite[Theorem~5.3.15]{kinglouis}. This amounts to showing that the structure sheaf $\scrO_{\M_\F(\Gamma)}$ evaluates to $K$ on the generic fiber and to $R$ globally. But $\scrO_{\M_\F(\Gamma)}(\M_\F(\Gamma)_\eta) = \scrO_{\F(V)}(\F(V))$, which equals $K$ by \cite[Corollary~3.3.21]{kinglouis}, whence the former claim. For the latter claim, note that $\scrO_{\M_\F(\Gamma)}(\M_\F(\Gamma))$ is a finitely generated $R$-module by \cite[Theorem~5.3.2\,(a)]{kinglouis}, and it is contained  in $\scrO_{\F(V)}(\F(V)) = K$, since $\M_\F(\Gamma)$ is integral, and therefore equal to $R$, since $R$ is integrally closed.

Equidimensionality of the special fiber follows from \cite[Corollaire~14.2.2]{EGAIV3}.
\end{Proof}

\begin{Remark}\label{rem:reduced}
  It has been proved in \cite[Theorem~2.3]{mustafinvarieties} that for $\F = \bbP$, the special fiber of a Mustafin degeneration  is always reduced. The methods used there do not readily apply to arbitrary flag varieties. However, since any Mustafin degeneration is a degeneration of the diagonal in $\F(V)^n$, up to a change of coordinates, the special fiber is generically reduced as long as this diagonal is \emph{multiplicity-free} as defined, e.g., in the introduction of \cite{brion-multfree}. This is a condition on $\F$, $d = \dim V$, and $n = |\Gamma|$. \cref{tab:multiplicity-free} lists some instances where this is the case; one example where it is \emph{not} the case is $\F = \Gr_2$, $d = 4$ and $n = 4$.

  \begin{table}
    \centering
    \begin{tabular}{lll}
      \toprule
      Flag type $\F$ & $d = \dim V$ & $n = |\Gamma|$ \tabularnewline
      \midrule
      arbitrary & arbitrary & $2$ \tabularnewline
      $\bbP$ & arbitrary & arbitrary \tabularnewline
      $(1 < 2)$ & $3$ & arbitrary \tabularnewline
      \bottomrule
    \end{tabular}
    \caption{Some cases where the diagonal in $\F(V)^n$ is multiplicity-free.}
    \label{tab:multiplicity-free}
  \end{table}
\end{Remark}


%% file: combinatorics.tex
\section{Components of the Special Fiber}\label{sec:combinatorics}

If $\Gamma \subseteq \Gamma'$ are finite subsets of $\frB_d^0$, then the following commutative diagram gives a natural surjective morphism $\M_\F(\Gamma') \twoheadrightarrow \M_\F(\Gamma)$:

\[ \xymatrix{& \M_\F(\Gamma') \ar@{^(->}[r] \ar@{-->>}[dd] & \prod_{[L] \in \Gamma'} \F(L) \ar@{->>}[dd] \\  
  \F(V) \ar@{^(->}[ur] \ar@{^(->}[dr] \\
  & \M_\F(\Gamma) \ar@{^(->}[r] & \prod_{[L] \in \Gamma} \F(L) } \]

The existence of the dashed morphism follows from the fact that under the projection $\prod_{[L] \in \Gamma'} \F(L) \twoheadrightarrow \prod_{[L] \in \Gamma} \F(L)$, the image of $\F(V)$ in the former product is mapped to the image of $\F(V)$ in the latter product, and hence, by continuity, its closure in $\prod_{[L] \in \Gamma'} \F(L)$ is mapped to its closure in $\prod_{[L] \in \Gamma} \F(L)$. Surjectivity follows from the fact that the morphism is dominant, since its image contains the generic fiber, and projective, hence closed.

From now on, denote by $X_s$ the special fiber of a scheme $X$ over $\Spec R$.

\begin{Lemma}\label{lem:birational}
  Let $\Gamma \subseteq \Gamma'$ be finite subsets of $\frB_d^0$. For each irreducible component $C$ of $\M_\F(\Gamma)_s$, there is a unique irreducible component $C'$ of $\M_\F(\Gamma')_s$ that maps onto $C$ via the natural projection $\M_\F(\Gamma') \twoheadrightarrow \M_\F(\Gamma)$. Furthermore, the map $C' \twoheadrightarrow C$ is birational.
\end{Lemma}

\begin{Proof}
  First, let us note that the special fiber of $\M_\F(\Gamma)$ has codimension $1$, because it corresponds to the prime ideal $\langle \pi \rangle$ for a uniformizer $\pi$ of $R$, and the only prime ideal contained in $\langle \pi \rangle$ is $\langle 0 \rangle$. By \cite[Corollary~4.4.3(b)]{kinglouis}, there is a non-empty open subset $V \subseteq \M_\F(\Gamma)$ such that the natural projection $\M_\F(\Gamma') \twoheadrightarrow \M_\F(\Gamma)$ is an isomorphism over $V$, and $\M_\F(\Gamma) \setminus V$ has codimension $\geq 2$. Therefore, $\M_\F(\Gamma) \setminus V$ has codimension $1$ in the special fiber. The claim follows from this.
\end{Proof}

\begin{Definition}
  If $\Gamma$ is a finite subset of $\frB_d^0$, then, for any $[L] \in \Gamma$, we call the unique component of $\M_\F(\Gamma)_s$ that maps birationally to $\M_\F(L)_s = \F(L)_s$ the \emph{primary component of $\M_\F(\Gamma)_s$ corresponding to $[L]$ (or $L$)}\index{component!primary}, or simply \emph{$L$-primary}. An irreducible component $C$ of $\M_\F(\Gamma)_s$ not mapping birationally to any $\M_\F(L)_s$ with $[L] \in \Gamma$ is called a \emph{secondary component} if there is a finite subset $\Gamma'$ of $\frB_d^0$ containing $\Gamma$, such that some primary component of $\M_\F(\Gamma')_s$ is mapped birationally to $C$. Any other irreducible component of $\M_\F(\Gamma)_s$ is called \emph{tertiary}.\index{component!tertiary}
\end{Definition}

In \cite[Lemma~5.8]{mustafinvarieties}, we showed that for $\F=\bbP$,  all components are primary or secondary. In what follows, we will show that tertiary components can appear for more general flag types.

\begin{Lemma}[Flag Projection Lemma]\label{lem:flag-projection}\index{flag projection}
  Let $\F' \subseteq \F$ be flag types and $\Gamma \subset \frB_d^0$ a finite subset. There is a natural surjective morphism $\M_\F(\Gamma) \twoheadrightarrow \M_{\F'}(\Gamma)$, under which, for any $[L] \in \Gamma$, the primary component of $\M_\F(\Gamma)_s$ corresponding to $L$ is mapped onto the primary component of $\M_{\F'}(\Gamma)_s$ corresponding to $L$.
\end{Lemma}

\begin{Proof}
  Since under the natural projection $\prod_{[L] \in \Gamma} \F(L) \twoheadrightarrow \prod_{[L] \in \Gamma} \F'(L)$, the diagonal of the generic fiber on the left is mapped onto the diagonal of the generic fiber on the right, we obtain the desired natural surjective morphism  $\M_\F(\Gamma) \twoheadrightarrow \M_{\F'}(\Gamma)$.

  Now consider the following cube of morphisms:
   \[\xymatrix{
      & C \ar@{~}[dd] |\hole  \ar@{^(->}[dl] \ar @ {.>>} [rr] & & C' \ar@{~}[dd]  \ar@{^(->}[dl] \\
      \M_\F(\Gamma)_s \ar@{^(->}[dd] \ar@{-->} [rr] && \M_{\F'}(\Gamma)_s \ar@{^(->}[dd] \\
      & \F(L)_s \ar@{->>} [rr] |(.49)\hole && \F'(L)_s \\
      \prod\limits_{[L] \in \Gamma} \F(L)_s \ar@{->>}[rr] \ar@{->>}[ur] && \prod\limits_{[L] \in \Gamma} \F'(L)_s \ar@{->>}[ur]
   }\]

   $C$ and $C'$ are the primary components corresponding to $L$ in $\M_\F(\Gamma)_s$ and $\M_{\F'}(\Gamma)_s$, respectively, the squiggly lines $\xymatrix{\ar@{~}[r]&}$ indicate birationality, and the left and the right face of the cube are commutative by definition; the bottom face consists of natural projections and is trivially commutative; the dashed arrow $\xymatrix{\ar@{-->}[r]&}$ (making the front face commutative) is the special fiber of the natural projection from above; and from this follows the existence of the dotted arrow $\xymatrix{\ar@{.>>}[r]&}$ making everything commutative, which is dominant because the other arrows on the back face are, and therefore is surjective.
\end{Proof}

We refer to this type of projection as \emph{flag projection}, as opposed to the projection defined at the beginning of this section stemming from an inclusion of vertex sets $\Gamma \subseteq \Gamma'$.

\begin{Lemma}\label{lem:unique-primaries}
  Let $\Gamma \subset \frB_d^0$ be finite, $[L_1], [L_2] \in \Gamma$ two distinct vertices, and $C_1, C_2$ the corresponding primary components in $\M_\F(\Gamma)_s$. Then $C_1 \neq C_2$.
\end{Lemma}

\begin{Proof}
  Both projections $\M_\F(\Gamma) \twoheadrightarrow \M_\F(L_1)$ and $\M_\F(\Gamma) \twoheadrightarrow \M_\F(L_2)$ factor through $\M_\F(\Gamma) \twoheadrightarrow \M_\F(L_1, L_2)$, so it suffices to show that the images of $C_1$ and $C_2$ in $\M_\F(L_1, L_2)$ are different. So assume $\Gamma = \{ [L_1], [L_2] \}$.
  Now if $\F = \Gr_k$ for some $k$, an explicit calculation using the equations given in \cref{sec:definitions} shows that the (unique) components mapping birationally to $\M_{\Gr_k}(L_1)$ and $\M_{\Gr_k}(L_2)$, respectively, are distinct.
 But if, for general $\F$, we had $C_1 = C_2$, 
then by the preceeding Flag Projection Lemma, the same would be true in $\M_{\Gr_k}(L_1, L_2)$ for any $\Gr_k \subset \F$, which we just ruled out.
\end{Proof}

\begin{Lemma}
  Let $\Gamma \subset \frB_d^0$ be finite, $C \subset \M_\F(\Gamma)_s$ a secondary irreducible component, and let $\Gamma'$ and $\Gamma''$ be two finite subsets of $\frB_d^0$ such that $C$ becomes primary in both $\M_\F(\Gamma')_s$ and $\M_\F(\Gamma'')_s$. Explicitly, this means that the unique irreducible components $C'$ and $C''$  of $\M_\F(\Gamma')_s$ and $\M_\F(\Gamma'')_s$, respectively, which map  birationally to $C$, are $L'$-primary and $L''$-primary, respectively, for some  $[L'] \in \Gamma'$ and $[L''] \in \Gamma''$.
Then we have $[L'] = [L'']$.
\end{Lemma}

\begin{Proof}
  Let $\tilde \Gamma = \Gamma' \cup \Gamma''$, and let $\tilde C'$ and $\tilde C''$ be the components of $\M_\F(\tilde \Gamma)_s$ mapping to $C'$ and $C''$ under the respective projections. We thus have the following diagram, squiggly lines $\xymatrix{\ar@{~}[r]&}$ indicating birationality:

\[ \xymatrix{
  \tilde C'\; \ar@{^(->}[rr] \ar@{~}[ddd] && \M_\F(\tilde \Gamma) \ar@{->>}[dl] \ar@{->>}[dr] \ar@{->>}[dd] && \; \tilde C'' \ar@{_(->}[ll] \ar@{~}[ddd] \\
 & \M_\F(\Gamma') \ar@{->>}[d] \ar@{->>}[dr] & & \M_\F(\Gamma'') \ar@{->>}[d] \ar@{->>}[dl] &  \\
 & \M_\F(L') & \M_\F(\Gamma) &  \M_\F(L'') &   \\
 \bigl.C'\bigr. \ar@{~}[rr]  \ar@/^/@{_(->}[uur] \ar@{ ~}[ur] &*{}& \; \bigl.C\bigr. \; \ar@<-1.5pt>@{^ (->}[u] && \bigl.C''\bigr. \ar@/_/@{^(-<}[luu] \ar@{~}[ul] \ar@{~}[ll] & 
}\] 

Now $\tilde C'$ and $\tilde C''$ both map birationally to $C$, so they are equal. But since they are primary components of $\M_\F(\tilde \Gamma)_s$ corresponding to $[L']$ and $[L'']$ respectively, $[L']$ and $[L'']$ must be equal.
\end{Proof}

This allows us to make the following definition:

\begin{Definition}
  Let $\Gamma \subset \frB_d^0$ be finite, und $C$ a secondary component of $\M_\F(\Gamma)_s$. We call $C$ \emph{the secondary component corresponding to $[L]$ (or $L$)}, or simply \emph{$L$-secondary}, if $[L]$ is the vertex in $\frB_d^0$ uniquely determined by the property that the $L$-primary component of $\M_\F(\Gamma \cup \{[L]\})$ maps birationally to $C$.
\end{Definition}

We have seen in \cref{lem:flag-projection} that $L$-primarity is preserved under flag projections. The same is true for $L$-secondarity, but only under an additional precondition:

\begin{Proposition}\label{prop:sec-comps-under-flag-projection}
  Let $\Gamma \subset \frB_d^0$ be finite, $\F' \subseteq \F$ flag types. Let $C$ be a secondary component of $\M_\F(\Gamma)_s$ corresponding to $[L] \in \frB_d^0$. Then $C$ is mapped onto an irreducible component $C'$ of $\M_{\F'}(\Gamma)_s$ by the flag projection if and only if there is a secondary component in $\M_{\F'}(\Gamma)_s$ corresponding to $[L]$. In this case, $C'$ is mapped to this $L$-secondary  component.
\end{Proposition}

\begin{Proof}
  Let $\tilde C$ be the $L$-primary component of $\M_\F(\Gamma \cup \{L\})_s$, which is mapped birationally to $C$ under the projection $\M_\F(\Gamma \cup \{L\}) \twoheadrightarrow \M_\F(\Gamma)$. By \cref{lem:flag-projection}, $\tilde C$ is mapped onto the $L$-primary component $\tilde C'$ of $\M_{\F'}(\Gamma \cup \{L\})_s$. By commutativity of the diagram
\[ \xymatrix{
  \M_\F(\Gamma \cup \{L\}) \ar@{->>}[r] \ar@{->>}[d] & \M_{\F'}(\Gamma \cup \{L\}) \ar@{->>}[d] \\
  \M_\F(\Gamma) \ar@{->>}[r] & \M_{\F'}(\Gamma)\text{,}
}\]
$\tilde C'$ is mapped onto a component $C'$ of $\M_{\F'}(\Gamma)_s$ if and only if $C$ is. Now if there is an $L$-secondary component in $\M_{\F'}(\Gamma)_s$, it is unique and $\tilde C'$ must be mapped onto it. Conversely, if $\tilde C'$ is mapped onto a component $C'$, it must be birationally by \cref{lem:birational}, making $C'$ an $L$-secondary component.
\end{Proof}

\begin{Definition}
  An irreducible component $C$ in the special fiber of a Mustafin degeneration $\M_\F(\Gamma)$ is called \emph{mixed} if there are flag types $\F_1, \F_2 \subseteq \F$ and vertices $[L_1], [L_2] \in \frB_d^0$, $[L_1] \neq [L_2]$, such that $C$ is mapped onto a (primary or secondary) component corresponding to $L_i$ under the flag projection to $\M_{\F_i}(\Gamma)$, for $i \in \{1,2\}$.
\end{Definition}

\cref{lem:flag-projection,prop:sec-comps-under-flag-projection} together imply that every mixed component is tertiary.

  The condition that $C$ is mapped onto an irreducible component in \cref{prop:sec-comps-under-flag-projection} is not automatic. In fact, the following example shows that there might be no secondary components at all in $\M_{\F'}(\Gamma)_s$ for $C$ to be mapped to!

\begin{Example}[Secondary and mixed components]\label{ex:mixed-component}
  Take $R = \bbQ[t]_{(t)}$ with field of fractions $K = \bbQ(t)$ and residue field $\bbQ$, $d = 3$, and $\F = ( 1 < 2 )$. For a basis $\{e_1, e_2, e_3\}$ of $V$, let 
  \begin{align*}
    L_1 &= Re_1 + Re_2 + Re_3\text{,} \\
    L_2 &= R\pi e_1 + Re_2 + Re_3\text{,} \\
    L_3 &= Re_1 + Re_2 + R\pi e_3\text{,} \\
  \end{align*}
  and $\Gamma = \{[L_1], [L_2], [L_3]\}$. By an explicit calculation, using the equations from \cref{sec:definitions}, we find that $\M_\F(\Gamma)_s$ possesses one secondary component $C$ corresponding to $L_4 = R\pi e_1 + Re_2 + R\pi e_3$, four mixed components and (of course) three primary components.

  In this example, the non-primary components are all toric threefolds: the secondary component $C$ is isomorphic to $(\bbP^1_\bbQ)^3$, two of the mixed components are isomorphic to $\bbP^1_\bbQ \times \widetilde \bbP^2_\bbQ$, where $\widetilde\bbP^2_\bbQ$ is the blow-up of of $\bbP^2_\bbQ$ in a point, and the remaining two mixed components are mutually isomorphic non-singular toric threefolds corresponding to a fan whose intersection with the unit sphere results in a pentagonal bipyramid (with $f$-vector $(7,15,10)$). 
Note, however, that the whole Mustafin degeneration $\M_\F(\Gamma)$ is not toric since the primary components are not toric.

We can depict the configuration $\Gamma \cup \{[L_4]\}$ as follows in the apartment corresponding to the basis $\{e_1, e_2, e_3\}$, which is the subcomplex of $\frB_d$ spanned by all vertices of the form $[\sum_iR\pi^{a_i}e_i]$: 
  \begin{center}
    \begin{tikzpicture}
      \draw[step=1cm,gray,very thin] (-2.4,-.4) grid (2.4,1.9);
      \begin{scope}
        \clip (-2.4,-.4) rectangle (2.4,1.9);
        \foreach \x in {-3,...,5} \draw[gray,very thin,rotate=45] (-10,\x*0.707106781) -- (10,\x*0.707106781); 
      \end{scope}
      \filldraw (0,0) circle (3pt) node[anchor=north west] {$L_1$};
      \filldraw (-1,0) circle (3pt) node[anchor=north east] {$L_2$};
      \filldraw (1,1) circle (3pt) node[anchor=south west] {$L_3$};
      \filldraw[fill=white] (0,1) circle (3pt) node[anchor=south east] {$L_4$};
    \end{tikzpicture}
  \end{center}

By \cite[Theorem~2.10]{mustafinvarieties}, $\M_{\bbP^\vee}(\Gamma)_s$ does not contain secondary components, since $\Gamma$ is convex in the sense of \cite[Theorem~2.10]{mustafinvarieties}. Therefore, $C$ must be mapped to a proper irreducible closed subset of an irreducible component of $\M_{\bbP^\vee}(\Gamma)_s$ under the flag projection $\M_\F(\Gamma) \twoheadrightarrow \M_{\bbP^\vee}(\Gamma)$. However, for the projective space $\bbP$, where the dual notion of convexity has to be applied, $\Gamma$ is not convex, but has $\Gamma \cup \{[L_4]\}$ as its convex hull, so $\M_{\bbP}(\Gamma)_s$ contains a secondary component correspondig to $L_4$, and $C$ is mapped onto this component.

The mixed components in this example are identified by the fact that they map to components (primary or secondary) corresponding to different vertices in $\M_{\bbP}(\Gamma)_s$ and $\M_{\bbP^\vee}(\Gamma)_s$, respectively.
\end{Example}

\begin{Example}
  The example above shows another interesting phenomenon. If we remove $[L_1]$ from $\Gamma$ and consider $\Gamma' = \{[L_2], [L_3]\}$, rather than observing two secondary components corresponding to $L_1$ and $L_4$, as one might expect, we only get four mixed components (besides the obligatory primary components).
\end{Example}

We conclude this \lcnamecref{sec:combinatorics} by giving upper and lower bounds for the total numbers of components of a Mustafin degeneration.

\begin{Lemma}\label{lem:number-of-components}
  Let $\F = (k_1 < \dots < k_r)$ be a flag type and $\Gamma \subset \frB_d^0$ finite with $|\Gamma| = n$. Denote by $c$ the number of irreducible components in the special fiber of $\M_\F(\Gamma)$. Then we have the lower bound $c \geq n$. If $n=2$, we also have the upper bound
  \begin{align*} c &\leq \binom{d}{d-k_r, \; k_r-k_{r-1}, \; \ldots, \; k_2-k_1, \; k_1} \\
    &= \text { the number of Schubert cells in $\F(V)$. }
  \end{align*}
\end{Lemma}

\begin{Proof}
  $c \geq n$ holds because there is one primary component for each vertex in $\Gamma$, and primary components corresponding to different vertices cannot coincide by \cref{lem:unique-primaries}.

  In order to prove the upper bound for $n=2$, we consider the class of $\M_\F(\Gamma)_s$ in the Chow ring $\CH(\F(\kappa^d)^2)$, where $\kappa$ denotes the residue field of $R$. 
Up to a change of coordinates, $\M_\F(\Gamma)_K$ is embedded in $\F(V)^2$ as the diagonal, so their classes in the Chow ring $\CH(\F(V)^2)$ are the same. 
Since $\M_\F(\Gamma)_s$ is a specialization of $\M_\F(\Gamma)_K$, as in \cite[Section~20.3]{fulton}, they have the same class in $\CH(\F(\kappa^d)^2)$.
Now each irreducible component of $\M_\F(\Gamma)_s$ adds a term with non-negative coefficients, so there can at most be as many components as the sum of the coefficients in $[\Delta]$.
  Since $n = 2$, this sum is exactly the number of Schubert varieties by the duality property of Schubert classes. That this number equals the stated multinomial coefficient follows from the definition of the Schubert varieties.
\end{Proof}

\begin{Question}
  The upper bound in the $n=2$ case is sharp for $\F = \bbP$ by \cite[Proposition~4.6]{mustafinvarieties}. It is attained for vertices in \emph{general position} as defined before the cited proposition. We would be interested to know if this is also true for general flag types.
\end{Question}


%% file: main.bbl
\begin{thebibliography}{\foreignlanguage{russian}{Mus78}}

\bibitem[AB08]{abramenko-brown}
Peter Abramenko and Kenneth~S. Brown.
\newblock {\em Buildings. Theory and applications}.
\newblock Springer, New York, 2008.

\bibitem[Bri03]{brion-multfree}
Michel Brion.
\newblock Multiplicity-free subvarieties of flag varieties.
\newblock In: Luchezar~L. Avramov, editor, {\em Contemporary Mathematics},
  volume 331, pages 13--23. American Mathematical Society, 2003.

\bibitem[CHSW11]{mustafinvarieties}
Dustin Cartwright, Mathias H\"abich, Bernd Sturmfels and Annette Werner.
\newblock Mustafin varieties.
\newblock {\em Selecta Mathematica New Series}, \textbf{17}:757--793, 2011.

\bibitem[EGA IV 3]{EGAIV3}
Alexander Grothendieck and Jean Dieudonn{é}.
\newblock {\em {É}l{é}ments de G{é}om{é}trie Alg{é}brique: {IV.} {É}tude
  locale des schémas et des morphismes des schémas, troisième partie},
  volume~28 of {\em Publications Mathématiques de l'{I.H.É.S.}}
\newblock {I}nstitut des {H}autes {É}tudes Scientifiques, Bures-sur-Yvettes,
  1966.

\bibitem[Ful97]{youngtableaux}
William Fulton.
\newblock {\em Young tableaux: with applications to representation theory and
  geometry}.
\newblock London Mathematical Society student texts. Cambridge University
  Press, Cambridge, 1997.

\bibitem[Ful98]{fulton}
William Fulton.
\newblock {\em Intersection theory}.
\newblock Springer, Berlin, 1998.

\bibitem[Liu02]{kinglouis}
Qing Liu.
\newblock {\em Algebraic Geometry and Arithmetic Curves}.
\newblock Springer, Oxford, 2002.

\bibitem[Mum72]{mumford}
David Mumford.
\newblock An analytic construction of degenerating curves over complete local
  rings.
\newblock {\em Compositio Mathematica}, \textbf{24}:129--174, 1972.

\bibitem[\foreignlanguage{russian}{Mus78}]{mustafin}
\foreignlanguage{russian}{G. A. Mustafin}.
\newblock \foreignlanguage{russian}{Nearhimedova uniformizaciya}.
\newblock {\em \foreignlanguage{russian}{Matematicheskii0 Sbornik}},
  \textbf{105}, 1978.
\newblock (English translation in: G. A. Mustafin: \textit{Nonarchimedean
  Uniformization}. Math. USSR Sbornik \textbf{34}:2, 1978).

\end{thebibliography}
